# DIFFERENCE PROPHET INEQUALITIES FOR $[0,1]$-VALUED I.I.D. RANDOM VARIABLES WITH COST FOR OBSERVATIONS[1]

By Holger Kösters

*University of Münster*

Let $X_1, X_2, \ldots$ be a sequence of $[0,1]$-valued i.i.d. random variables, let $c \geq 0$ be a sampling cost for each observation and let $Y_i = X_i - ic$, $i = 1, 2, \ldots$. For $n = 1, 2, \ldots$, let $M(Y_1, \ldots, Y_n) = E(\max_{1 \leq i \leq n} Y_i)$ and $V(Y_1, \ldots, Y_n) = \sup_{\tau \in C^n} E(Y_\tau)$, where $C^n$ denotes the set of all stopping rules for $Y_1, \ldots, Y_n$. Sharp upper bounds for the difference $M(Y_1, \ldots, Y_n) - V(Y_1, \ldots, Y_n)$ are given under various restrictions on $c$ and $n$.

**1. Introduction.** In her interesting paper, Samuel-Cahn (1992) investigated "prophet inequalities" for $[0,1]$-valued i.i.d. random variables with cost for observations: Let $X_1, X_2, \ldots$ be i.i.d. random variables, $0 \leq X_i \leq 1$, and let $c \geq 0$ be a fixed sampling cost that is charged for each observation. Consider the sequence $Y_i = X_i - ic$, $i = 1, 2, \ldots$, and for $n = 1, 2, \ldots$, let $M(Y_1, \ldots, Y_n) = E(\max_{1 \leq i \leq n} Y_i)$ and $V(Y_1, \ldots, Y_n) = \sup_{\tau \in C^n} E(Y_\tau)$, where $C^n$ denotes the set of stopping rules for $Y_1, \ldots, Y_n$. Then $M(Y_1, \ldots, Y_n)$ and $V(Y_1, \ldots, Y_n)$ can be interpreted as the expected optimal return of a prophet and a statistician, respectively. For any real number $x$, let $[x]$ denote the largest integer strictly smaller than $x$. Samuel-Cahn (1992) stated her main result as follows:

Let $X_1, X_2, \ldots$ be i.i.d. random variables, $0 \leq X_1 \leq 1$.

(a) For $0 < c \leq 1$ fixed and all $n \geq 1$,
$$M(Y_1, \ldots, Y_n) - V(Y_1, \ldots, Y_n) \leq [1/c]c(1-c)^{[1/c]+1}.$$

(b) For $n \geq 1$ fixed and all $c \geq 0$,
$$M(Y_1, \ldots, Y_n) - V(Y_1, \ldots, Y_n) \leq (1 - 1/n)^{n+1}.$$

Received January 2003; revised October 2003.
[1]Supported by DFG Grant SCHM 677/8.

(c) For all $c \geq 0$ and all $n \geq 1$,
$$M(Y_1,\ldots,Y_n) - V(Y_1,\ldots,Y_n) < e^{-1}.$$

All bounds are the best possible.

Unfortunately, Harten (1996) detected a gap in the proof (the argument for the reduction to Bernoulli variables is incomplete) and even showed that the inequality in (a) fails to hold for $c > 1/2$ and that the inequality in (b) fails to hold for $n \geq 2$ (see also the remarks below). Moreover, because the original proof of result (c) was based on result (b), a new proof for this part became necessary as well.

Harten gave the correct upper bound for part (a) and provided a new proof for part (c) based on this result. With regard to part (b), Harten conjectured that the correct upper bound is $(n-1)(1-1/(n+1))^n/(n+1)$, but he proved this result only for the special case of Bernoulli variables (see Lemma 2.1). The purpose of the present note is to extend this result to arbitrary $[0,1]$-valued i.i.d. random variables. To summarize, we present the following complete result:

THEOREM. *Let $X_1, X_2, \ldots$ be i.i.d. random variables, $0 \leq X_1 \leq 1$.*

(a) *For $0 < c \leq 1$ fixed and all $n \geq 1$,*
$$M(Y_1,\ldots,Y_n) - V(Y_1,\ldots,Y_n) \leq \begin{cases} [1/c]c(1-c)^{[1/c]+1}, & \text{for } c \leq 1/2, \\ (1-c)/4, & \text{for } c \geq 1/2. \end{cases}$$

(b) *For $n \geq 1$ fixed and all $c \geq 0$,*
$$M(Y_1,\ldots,Y_n) - V(Y_1,\ldots,Y_n) \leq (n-1)(1-1/(n+1))^n/(n+1).$$

(c) *For all $c \geq 0$ and all $n \geq 1$,*
$$M(Y_1,\ldots,Y_n) - V(Y_1,\ldots,Y_n) < e^{-1}.$$

*All bounds are the best possible.*

The proof of part (a) can be found in the Ph.D. thesis of Harten (1996) or in Section 8(c) of Harten, Meyerthole and Schmitz (1997). It proceeds roughly as follows: In the case $EX_1 \leq c$, the problem can first be reduced to that for Bernoulli variables [i.e., random variables $X_i$ such that $P(X_i = 1) = p = 1 - P(X_i = 0)$ for some $0 \leq p \leq 1$] and then be solved by direct calculation. In the case $EX_1 > c$, the basic idea is to apply Theorem A from Jones (1990) to the sequence $V(Y_2,\ldots,Y_n), Y_2, \ldots, Y_n$ and to use some suitable estimates.

The proof of part (b), presented in the following section, is somewhat similar in that we also distinguish the cases $EX_1 \leq c$ and $EX_1 > c$. However, in the second case we use a completely different argument.

Part (c) follows easily from part (a) or part (b), since $e^{-1}$ is the supremum of the respective upper bounds.

REMARKS. (a) The following examples, taken from Harten (1996), show that the upper bounds are attained in the cases (a) and (b): In (a) take $n \geq [1/c] + 1$ and i.i.d. random variables $X_1, \ldots, X_n$ with $P(X_1 = 1) = c = 1 - P(X_1 = 0)$ for $c \leq 1/2$, and $P(X_1 = 1) = 1/2 = 1 - P(X_1 = 0)$ for $c \geq 1/2$. In (b) take $c = 1/(n+1)$ and i.i.d. random variables $X_1, \ldots, X_n$ with $P(X_1 = 1) = 1/(n+1) = 1 - P(X_1 = 0)$. We note without proof that, except for $c = 1$ in (a) and $n = 1$ in (b), the upper bounds are attained only in these cases.

(b) The inequalities in parts (a) and (c) remain true for infinite sequences of $[0, 1]$-valued i.i.d. random variables. See Harten (1996), Harten, Meyerthole and Schmitz (1997) or Saint-Mont (1999) for a more general presentation.

(c) It seems natural to look also for ratio prophet inequalities, that is, for inequalities of the type $M(Y_1, \ldots, Y_n)/V(Y_1, \ldots, Y_n) \leq C$. However, this ratio turns out to be unbounded for all $n > 1$ and all $0 < c < 1$, as already observed by Samuel-Cahn (1992).

(d) The case $c = 0$ was treated by Hill and Kertz (1982), who showed that $M(Y_1, \ldots, Y_n)/V(Y_1, \ldots, Y_n) \leq a_n$ and $M(Y_1, \ldots, Y_n) - V(Y_1, \ldots, Y_n) \leq b_n$ for all $n = 2, 3, \ldots$, with certain constants $1.1 < a_n < 1.6$ and $0 < b_n < 1/4$. These results are markedly different from those for the case $c > 0$ and seemingly cannot be obtained from them. Quite on the contrary, the second inequality is a key ingredient in our proof of part (b).

**2. Proving part (b).** Throughout this section, assume that $n \geq 2$ (otherwise the assertion is trivial), and let $D(Y_1, \ldots, Y_n) := M(Y_1, \ldots, Y_n) - V(Y_1, \ldots, Y_n)$ and $d_n := (n-1)(1 - 1/(n+1))^n/(n+1)$. It remains to prove that

$$D(Y_1, \ldots, Y_n) \leq d_n. \tag{1}$$

The examples given at the end of the Introduction then show that this bound is also the best possible.

In two special cases, the Bernoulli case and the zero-cost case, (1) can easily be deduced from existing results:

LEMMA 2.1. *For all $c \geq 0$ and all i.i.d. random variables $X_1, \ldots, X_n$ with $P(X_1 = 1) = 1 - P(X_1 = 0)$, $D(Y_1, \ldots, Y_n) \leq d_n$.*

PROOF. For $c = 0$ and $c \geq 1$, $D(Y_1, \ldots, Y_n) = 0 < d_n$, since the statistician can secure the same return as the prophet by using the stopping rules $\tau = \inf\{i | X_i = 1 \text{ or } i = n\}$ and $\tau \equiv 1$, respectively. For $0 < c < 1$, the result is obtained by direct calculation, and can be found in Harten [(1996), pages 142 and 143] or Harten, Meyerthole and Schmitz [(1997), pages 194–196]. □

LEMMA 2.2. *Suppose that $c = 0$. Then for all i.i.d. random variables $X_1, \ldots, X_n$ taking values in $[0,1]$, $D(Y_1, \ldots, Y_n) < d_n$.*

PROOF. For $c = 0$, Theorem B in Hill and Kertz (1982) states that

$$D(Y_1, \ldots, Y_n) = D(X_1, \ldots, X_n) \leq b_n$$

for certain constants $0 < b_n < 1/4$ [see Hill and Kertz (1982) for definitions]. Hence it remains to show that $b_n < d_n$: For $n \geq 5$ this follows from $d_n > 1/4$. For $n \leq 4$, we can go back to the definitions to verify the inequalities $b_2 \cong 0.063 < 0.148 \cong d_2$, $b_3 \cong 0.077 < 0.211 \cong d_3$ and $b_4 \cong 0.085 < 0.246 \cong d_4$ [see also Example 3.9 in Hill and Kertz (1982)]. □

We now show that (1) remains true for any $c \geq 0$ and any $[0,1]$-valued i.i.d. random values $X_1, \ldots, X_n$ by relating the general case to the above-mentioned special cases.

In doing so, we use that, from Theorem 3.2 in Chow, Robbins and Siegmund (1971),

$$(2) \qquad V(Y_1, \ldots, Y_i) = E(\max\{X_1, V(Y_1, \ldots, Y_{i-1})\}) - c$$

for all $i > 1$. Furthermore, we use the usual balayage technique [the reduction to distributions with maximum variance; see Section 2 in Hill and Kertz (1982)]. For any integrable random variable $Y$ and any $-\infty < a < b < \infty$, let $Y_a^b$ denote a random variable with $Y_a^b = Y$ for $Y \notin [a,b]$, $= a$ with probability $(b-a)^{-1} \int_{Y \in [a,b]} (b - Y) \, dP$ and $= b$ with probability $(b-a)^{-1} \int_{Y \in [a,b]} (Y - a) \, dP$. Then $EY_a^b = EY$ and if $X$ is any random variable independent of $Y$ and $Y_a^b$,

$$(3) \qquad E(\max\{X, Y\}) \leq E(\max\{X, Y_a^b\}).$$

PROOF OF PART (b). For $c = 0$, the assertion has just been established.

For $c > 0$ and $EX_1 \leq c$, we follow Harten [(1996), Proposition 12.4] and reduce the problem directly to that for Bernoulli variables. Let $\tilde{X}_1, \ldots, \tilde{X}_n$ be i.i.d. random variables, where $\tilde{X}_1 := (X_1)_0^1$ is a 0–1 balayage of $X_1$ and independent of $X_1, \ldots, X_n$, and let $\tilde{Y}_i = \tilde{X}_i - ic$, $i = 1, \ldots, n$. Then it follows from (2) and (3) that $V(\tilde{Y}_1, \ldots, \tilde{Y}_i) = V(Y_1, \ldots, Y_i)$ for all $i = 1, \ldots, n$ and $M(\tilde{Y}_1, \ldots, \tilde{Y}_n) \geq M(Y_1, \ldots, Y_n)$. Therefore, $D(Y_1, \ldots, Y_n) \leq D(\tilde{Y}_1, \ldots, \tilde{Y}_n)$, but now $\tilde{X}_1, \ldots, \tilde{X}_n$ are Bernoulli variables, so referring to Lemma 2.1 yields (1).

For $c > 0$ and $EX_1 > c$, a direct reduction to Bernoulli variables does not seem possible, so we have to take a different approach. For better clarity, we split the proof into several steps:

First of all we show that we may restrict ourselves to i.i.d. $[0,1]$-valued random variables $X_1, \ldots, X_n$ such that

$$(4) \qquad \begin{aligned} x_* &:= \inf\{x \in \mathbb{R} | P(X_1 \leq x) > 0\} = 0, \\ x^* &:= \sup\{x \in \mathbb{R} | P(X_1 \leq x) < 1\} = 1, \end{aligned}$$

that is, the length of the interval $[0,1]$ is fully exhausted:

LEMMA 2.3. *Let $X_1, \ldots, X_n$ be i.i.d. $[0,1]$-valued random variables and let $c > 0$. Then there exist i.i.d. $[0,1]$-valued random variables $\tilde{X}_1, \ldots, \tilde{X}_n$ with*

$$\tilde{x}_* := \inf\{\tilde{x} \in \mathbb{R} | P(\tilde{X}_1 \leq \tilde{x}) > 0\} = 0,$$
$$\tilde{x}^* := \sup\{\tilde{x} \in \mathbb{R} | P(\tilde{X}_1 \leq \tilde{x}) < 1\} = 1$$

*and $\tilde{c} > 0$ such that $D(Y_1, \ldots, Y_n) \leq D(\tilde{Y}_1, \ldots, \tilde{Y}_n)$. Here let $\tilde{Y}_i$ be defined by $\tilde{Y}_i := \tilde{X}_i - i\tilde{c}$, $i = 1, \ldots, n$.*

PROOF. Let $x_*, x^*$ be defined as in (4).

For $x_* = x^*$, we have $M(Y_1, \ldots, Y_n) - V(Y_1, \ldots, Y_n) = 0$, and choosing i.i.d. random variables $\tilde{X}_1, \ldots, \tilde{X}_n$ with $P(\tilde{X}_i = 1) = 1/2 = 1 - P(\tilde{X}_i = 0)$ and $\tilde{c} := c$ yields the assertion.

For $x_* < x^*$, consider the random variables $\tilde{X}_1, \ldots, \tilde{X}_n$ defined by $\tilde{X}_i := (X_i - x_*)/(x^* - x_*)$. These are obviously a.s. $[0,1]$-valued i.i.d. random variables with $\tilde{x}_* = 0$ and $\tilde{x}^* = 1$, and setting $\tilde{c} := c/(x^* - x_*)$ $[> 0]$, we also have

$$M(\tilde{Y}_1, \ldots, \tilde{Y}_n) - V(\tilde{Y}_1, \ldots, \tilde{Y}_n)$$
$$= (M(Y_1, \ldots, Y_n) - x_*)/(x^* - x_*) - (V(Y_1, \ldots, Y_n) - x_*)/(x^* - x_*)$$
$$\geq M(Y_1, \ldots, Y_n) - V(Y_1, \ldots, Y_n).$$

Thus, after a modification on a null set if necessary, the random variables $\tilde{X}_i$ have the desired properties. □

Note that in the preceding reduction step, we possibly get from the case $c > 0$, $EX_1 > c$ to the case $\tilde{c} > 0$, $E\tilde{X}_1 \leq \tilde{c}$. This being supposed, the first part of the proof yields the assertion. Hence it remains to consider the case $\tilde{c} > 0$, $E\tilde{X}_1 > \tilde{c}$.

Next we follow Hill and Kertz [(1982), Lemma 2.4] and show that we may restrict attention to special discrete distributions. Indeed, it suffices to consider i.i.d. $[0,1]$-valued random variables $X_1, \ldots, X_n$ such that

(5) $\quad\quad P(X_1 \in \{1, v_{n-1}, v_{n-2}, \ldots, v_2, v_1, 0\}) = 1,$
$\quad\quad\quad\, P(X_1 = 1) > 0, \quad P(X_1 = 0) > 0,$

where $v_i := V(Y_1, \ldots, Y_i)$, $i = 1, \ldots, n-1$.

LEMMA 2.4. *Let $X_1, \ldots, X_n$ be i.i.d. $[0,1]$-valued random variables satisfying condition (4), let $c > 0$ and suppose that $EX_1 > c$. Let $v_i := V(Y_1, \ldots, Y_i)$,* $i = 1, \ldots, n - 1$. Then $1 > v_{n-1} > v_{n-2} > \cdots > v_2 > v_1 > 0$ and there exist i.i.d. $[0,1]$-valued random variables $\tilde{X}_1, \ldots, \tilde{X}_n$ with

$$P(\tilde{X}_1 \in \{1, v_{n-1}, v_{n-2}, \ldots, v_2, v_1, 0\}) = 1,$$
$$P(\tilde{X}_1 = 1) > 0, \qquad P(\tilde{X}_1 = 0) > 0$$

and

$$V(\tilde{Y}_1, \ldots, \tilde{Y}_i) = v_i = V(Y_1, \ldots, Y_i) \qquad \text{for all } i = 1, \ldots, n-1$$

such that $D(Y_1, \ldots, Y_n) \leq D(\tilde{Y}_1, \ldots, \tilde{Y}_n)$. Here let $\tilde{Y}_i$ be defined by $\tilde{Y}_i := \tilde{X}_i - ic$, $i = 1, \ldots, n$.

PROOF. We begin by proving the inequality $1 > v_{n-1} > v_{n-2} > \cdots > v_2 > v_1 > 0$. Clearly, $v_i \leq 1 - c < 1$ for all $i = 1, \ldots, n-1$. Furthermore, $v_1 = EX_1 - c > 0$ by assumption. Moreover, if $v_i > v_{i-1}$ holds for some $i \in \{1, \ldots, n-2\}$ (where $v_0 := 0$), (2) yields $v_{i+1} = E(\max\{X_1, v_i\}) - c > E(\max\{X_1, v_{i-1}\}) - c = v_i$, where the strict inequality follows from the assumption $x_* = 0$.

Now, using the balayage technique, choose i.i.d. random variables $\tilde{X}_1, \ldots, \tilde{X}_n$ with the same distribution as $((\cdots((X_1)_0^{v_1})_{v_1}^{v_2} \cdots)_{v_{n-2}}^{v_{n-1}})_{v_{n-1}}^{1}$. Then it is obvious that

$$P(\tilde{X}_1 \in \{1, v_{n-1}, v_{n-2}, \ldots, v_2, v_1, 0\}) = 1,$$
$$P(\tilde{X}_1 = 1) > 0, \qquad P(\tilde{X}_1 = 0) > 0.$$

(For the inequalities, we need the assumptions $x^* = 1$ and $x_* = 0$.) Furthermore, the same argument as in the proof of Lemma 2.4 in Hill and Kertz (1982) shows that $V(\tilde{Y}_1, \ldots, \tilde{Y}_i) = V(Y_1, \ldots, Y_i)$ for all $i = 1, \ldots, n$ and $M(\tilde{Y}_1, \ldots, \tilde{Y}_n) \geq M(Y_1, \ldots, Y_n)$, whence $D(Y_1, \ldots, Y_n) \leq D(\tilde{Y}_1, \ldots, \tilde{Y}_n)$. □

Note that in this reduction step, passing from the $X_i$ to the $\tilde{X}_i$ leaves the expectation unchanged, so that we stay in the case $c > 0$, $EX_1 > c$.

Whereas in the i.i.d. case without observation costs, the reduction to special discrete distributions leads to a tractable formula for computing $M(Y_1, \ldots, Y_n)$ with the aid of the distribution function of $X_1$ [see Hill and Kertz (1982), Lemma 2.5], such a procedure does not seem to be possible in the i.i.d. case with observation costs. The reason for this is that we do not know enough about the order relationships between the $v_i - hc$, $i = 1, \ldots, n-1$, $h = 1, \ldots, n$.

To circumvent this problem, we embed the random variables $Y_1, \ldots, Y_n$ into a whole family of random variables $Y_1(\beta), \ldots, Y_n(\beta)$ in such a way that (i) we can easily bound the difference $D(Y_1(\beta), \ldots, Y_n(\beta))$ from above for two special values of the parameter $\beta$ and (ii) the resulting bounds lead to an upper bound for the original difference $D(Y_1, \ldots, Y_n)$.

LEMMA 2.5. *Let $X_1, \ldots, X_n$ be i.i.d. $[0,1]$-valued random variables and let $c > 0$ such that the conditions $EX_1 > c$ and (5) are satisfied. Then*

$$c' := c - P(X_1 = 1) < 0.$$

PROOF. Suppose by way of contradiction that $c' \geq 0$, that is, $P(X_1 = 1) \leq c$. Since with probability 1, $X_1$ takes on the values $0 < v_1 < \cdots < v_{n-1} < 1$ only, we obtain

$$\begin{aligned} v_{n-1} &= E(\max\{X_1, v_{n-2}\}) - c \\ &= v_{n-2} \cdot P(X_1 < v_{n-1}) + v_{n-1} \cdot P(X_1 = v_{n-1}) + P(X_1 = 1) - c \\ &\leq v_{n-1} \cdot P(X_1 < v_{n-1}) + v_{n-1} \cdot P(X_1 = v_{n-1}) < v_{n-1}, \end{aligned}$$

where the last inequality follows from $P(X_1 = 1) > 0$. This contradiction proves the lemma. □

CONSTRUCTION 2.6. Let $X_1, \ldots, X_n$ be i.i.d. $[0,1]$-valued random variables, let $c > 0$ such that the conditions $EX_1 > c$ and (5) are satisfied, and let $\beta^* := -P(X_1 = 1)/c'$ ($> 0$) (with $c'$ as in Lemma 2.5). We now construct a family of random variables $\{X_1(\beta), \ldots, X_n(\beta)\}_{\beta \in [0, \beta^*]}$ with corresponding sampling costs $\{c(\beta)\}_{\beta \in [0, \beta^*]}$. For all $\beta \in [0, \beta^*]$, let $v_i(\beta) := \beta \cdot v_i$, $i = 1, \ldots, n-1$,

$$X_h(\beta) := \begin{cases} 1, & \text{on } \{X_h = 1\}, \\ v_i(\beta), & \text{on } \{X_h = v_i\}, \ h = 1, \ldots, n, \\ 0, & \text{on } \{X_h = 0\}, \end{cases}$$

and $c(\beta) := \beta \cdot c' + P(X_1 = 1)$. Finally, let

$$Y_h(\beta) := X_h(\beta) - h \cdot c(\beta), \qquad h = 1, \ldots, n.$$

Note that $\beta^* = -P(X_1 = 1)/c'$ is the (uniquely determined) zero of the strictly decreasing function $\beta \mapsto \beta \cdot c' + P(X_1 = 1)$. Since $1 \cdot c' + P(X_1 = 1) = c > 0$, this implies $\beta^* > 1 > 0$. In particular, $\beta = 1$ is an admissible parameter. Furthermore, $c(\beta) \geq 0$ for all $\beta \in [0, \beta^*]$, with equality holding if and only if $\beta = \beta^*$.

LEMMA 2.7. *Given the situation of Construction 2.6, we have:*

(a) *For $\beta = 1$, the $X_i(\beta)$ and $c(\beta)$ coincide with the $X_i$ and $c$.*
(b) *For each $\beta \in (0, \beta^*)$, $X_1(\beta), \ldots, X_n(\beta)$ are i.i.d. $[0,1]$-valued random variables such that condition (5) is satisfied with respect to the sampling cost $c(\beta)$.*
(c) *For $\beta = 0$, $X_1(\beta), \ldots, X_n(\beta)$ are i.i.d. random variables with $P(X_1(\beta) = 1) = 1 - P(X_1(\beta) = 0)$ and $c(\beta) = P(X_1 = 1)$.*

(d) *For $\beta = \beta^*$, $X_1(\beta), \ldots, X_n(\beta)$ are i.i.d. $[0,1]$-valued random variables and $c(\beta) = 0$.*

PROOF. Statement (a) is obvious from the definitions.

To prove (b)–(d), let $\tilde{\beta} := \sup\{\beta \in [0, \beta^*] | v_{n-1}(\beta) \leq 1\}$. Since $v_{n-1}(1) = v_{n-1} < 1$, we have $\tilde{\beta} > 1$. Furthermore, it is obvious that

(6) $$0 \leq v_1(\beta) \leq \cdots \leq v_{n-1}(\beta) \leq 1$$

for all $\beta \in [0, \tilde{\beta}]$, where the first $n-1$ equalities hold exactly for $\beta = 0$ and the last equality holds at most for $\beta = \tilde{\beta}$.

We now show by induction on $h$ that

(7) $$V(Y_1(\beta), \ldots, Y_h(\beta)) = v_h(\beta) \qquad \text{for all } \beta \in [0, \tilde{\beta}]$$

for all $h = 1, \ldots, n-1$. First note that $X_1(\beta), \ldots, X_n(\beta)$ are again i.i.d. random variables and that with probability 1, $X_1(\beta) = \beta \cdot X_1 \cdot \mathbb{1}_{\{X_1 < 1\}} + \mathbb{1}_{\{X_1 = 1\}}$. Setting $v_0 := 0$, $v_0(\beta) := 0$ for $h = 1$ and using the inductive hypothesis $V(Y_1(\beta), \ldots, Y_{h-1}(\beta)) = v_{h-1}(\beta)$ for $h = 2, \ldots, n-1$, we therefore obtain

$$\begin{aligned}
V(Y_1(\beta), &\ldots, Y_h(\beta)) \\
&= E(\max\{X_1(\beta), v_{h-1}(\beta)\}) - c(\beta) \\
&= E(\max\{X_1(\beta), v_{h-1}(\beta)\} \cdot \mathbb{1}_{\{X_1 < 1\}}) + P(X_1 = 1) - \beta \cdot c' - P(X_1 = 1) \\
&= \beta \cdot (E(\max\{X_1, v_{h-1}\} \cdot \mathbb{1}_{\{X_1 < 1\}}) + P(X_1 = 1) - c' - P(X_1 = 1)) \\
&= \beta \cdot (E(\max\{X_1, v_{h-1}\}) - c) \\
&= \beta \cdot V(Y_1, \ldots, Y_h) = \beta \cdot v_h = v_h(\beta),
\end{aligned}$$

which proves (7).

We now show that $\tilde{\beta} = \beta^*$. By definition it is clear that $\tilde{\beta} \leq \beta^*$. Suppose by way of contradiction that $\tilde{\beta} < \beta^*$. Then for monotonicity and continuity reasons we have $0 < v_1(\tilde{\beta}) < \cdots < v_{n-1}(\tilde{\beta}) = 1$ and $c(\tilde{\beta}) > 0$ and therefore

$$1 = v_{n-1}(\tilde{\beta}) = V(Y_1(\tilde{\beta}), \ldots, Y_{n-1}(\tilde{\beta})) \leq 1 - c(\tilde{\beta}) < 1,$$

that is, a contradiction. Hence $\tilde{\beta} = \beta^*$.

It follows that (6) even holds for all $\beta \in [0, \beta^*]$. Since the other properties mentioned in (b)–(d) are obvious now, the proof is complete. □

LEMMA 2.8. *The function $\beta \mapsto V(\beta) := V(Y_1(\beta), \ldots, Y_n(\beta))$ is linear.*

PROOF. Similarly as in (7), we have $V(Y_1(\beta), \ldots, Y_n(\beta)) = \beta \cdot V(Y_1, \ldots, Y_n)$ for all $\beta \in [0, \beta^*]$. □

LEMMA 2.9. *The function $\beta \mapsto M(\beta) := M(Y_1(\beta), \ldots, Y_n(\beta))$ is convex.*

PROOF. For almost every $\omega \in \Omega$, the "path"

$$\beta \mapsto Y^*(\beta;\omega) := \max\{Y_1(\beta;\omega), Y_2(\beta;\omega), \ldots, Y_n(\beta;\omega)\}$$

is convex, since it is the maximum of the affine-linear functions $\beta \mapsto Y_i(\beta;\omega)$. Thus the function $\beta \mapsto M(\beta) = E(Y^*(\beta))$ is also convex. $\square$

CONTINUATION OF THE PROOF OF PART (b). We are now in a position to complete the proof of part (b). Let $X_1, X_2, \ldots$ and let $c$ be such that $c > 0$, $EX_1 > c$. By Lemmas 2.3 and 2.4, we may assume that condition (5) is satisfied (see also the remarks below Lemmas 2.3 and 2.4) and that Construction 2.6 is applicable. Then the function $\beta \mapsto D(\beta) := M(\beta) - V(\beta)$ is convex, since it is the difference of a convex and a linear function. Since a convex function defined on a compact interval always attains its maximum on the boundary of its domain, it follows that

(8) $$D(Y_1, \ldots, Y_n) = D(1) \leq \max\{D(0), D(\beta^*)\}.$$

Now, on the one hand, we have Bernoulli variables for $\beta = 0$ [Lemma 2.7(c)], which implies $D(0) \leq d_n$ by Lemma 2.1; on the other hand, we have $c(\beta) = 0$ for $\beta = \beta^*$ [Lemma 2.7(d)], which implies $D(\beta^*) < d_n$ by Lemma 2.2. Hence the maximum in (8) is bounded above by $d_n$, which proves (1). $\square$

**Acknowledgments.** I thank Professor N. Schmitz for the problem statement and for his friendly and constant support. Furthermore, I am grateful to the reviewers for their helpful comments on a previous version of this paper.

INSTITUT FÜR MATHEMATISCHE STATISTIK
UNIVERSITÄT MÜNSTER
EINSTEINSTRASSE 62
D-48149 MÜNSTER
GERMANY
E-MAIL: koestho@math.uni-muenster.de